\documentclass[12pt]{article}
\usepackage{amsfonts}
\usepackage{latexsym,amsmath,amssymb}
\usepackage{graphics,epstopdf}
\usepackage[pdftex]{graphicx}
\usepackage{color}
\usepackage{hyperref}
\hypersetup{colorlinks=false}
\usepackage[ top=1.5cm, bottom=1.5cm, outer=2cm, inner=2cm,
heightrounded, marginparwidth=2.5cm, marginparsep=2cm] {geometry}
\numberwithin{equation}{section}
\begin{document}

\bigskip

\bigskip

\begin{center}
{\Large \textbf{\ On a Kantorovich variant of $(p,q)$-Sz\'{a}sz-Mirakjan operators}}

\bigskip

\textbf{$^1$M. Mursaleen}{\footnote{\emph{Corresponding author}}}, \textbf{$^1$Khursheed J. Ansari} and  \textbf{$^2$Abylkassymova Elmira}

$^1$Department of\ Mathematics, Aligarh Muslim University, Aligarh-202002, India\\
$^2$chair "The theory and methods of teaching informatics",
Science-Pedagogical Faculty, M. Auezov South Kazakhstan State
University, Tauke Khan Avenue 5, Shymkent 160012, Kazakhstan
\\[0pt]
mursaleenm@gmail.com; ansari.jkhursheed@gmail.com; $smanchik_darmen@mail.ru$\\
\bigskip

\bigskip

\textbf{Abstract}
\end{center}

\parindent=8mm {\footnotesize {In the present paper we propose a Kantorovich variant of $(p,q)$-analogue of Sz\'{a}sz-Mirakjan operators.
We establish the moments of the operators with the help of a recurrence relation that we have derived
and then prove the basic convergence theorem. Next, the local approximation as well as
weighted approximation properties of these new operators in terms of modulus of continuity are studied.}}\newline


\\
\textbf{1. Introduction and Notations}\\

\parindent8mm Approximation theory has been an established field of mathematics in the past century.
The approximation of functions by positive linear operators is an important research topic
in general mathematics and it also provides powerful tools to application areas such as
computer-aided geometric design, numerical analysis, and solution of differential equations.

\parindent8mm During the last two decades, the applications of $q$-calculus emerged
as a new area in the field of approximation theory. The rapid development of
$q$-calculus has led to the discovery of various generalizations of
Bernstein polynomials involving $q$-integers. Several researchers introduced and studied
many positive linear operators based on $q$-integers, $q$-Bernstein basis, $q$-Beta basis, $q$-derivative and $q$-integrals etc.
Using $q$-integers, Lupa\c{s} \cite{lp} introduced the first $q$%
-Bernstein operators \cite{brn} and investigated its approximating and
shape-preserving properties. Another $q$-analogue of the Bernstein
polynomials is due to Phillips \cite{pl}. Since then several generalizations
of well-known positive linear operators based on $q$-integers have been
introduced and studied their approximation properties. Aral \cite{aral} and Aral and Gupta \cite{aral1} proposed
and studied some $q$-analogue of Sz\'{a}sz-Mirakjan operators \cite{sza}, defined on positive real axis.
Also Mahmudov \cite{mah} introduced a $q$-parametric Sz\'{a}sz-Mirakjan operators and studied
their convergence properties. Recently, approximation properties for King's type q-Bernstein–Kantorovich operators have been studied in \cite{mkk}.

\parindent8mm Very recently, Mursaleen et al applied $(p,q)$-calculus in approximation theory
and introduced the{\ $(p,q)$-analogue of Bernstein operators \cite{mka1,mkar}
and $(p,q)$-Bernstein-Stancu operators  \cite{mka2} and investigated their approximation properties. Also Acar \cite{acar} has introduced
$(p,q)$ parametric generalization of Sz\'{a}sz-Mirakjan operators. In the present work we have proposed
a Kantorovich variant of Sz\'{a}sz-Mirakjan operators and establish the moments of the operators with
the help of a recurrence relation that we have derived
and then prove the basic convergence theorem. Next, the local approximation as well as
weighted approximation properties of these new operators in terms of modulus of continuity are studied.

\parindent8mm The $(p,q)$-integer was introduced in order to generalize or unify several
forms of $q$-oscillator algebras well known in the earlier physics
literature related to the representation theory of single parameter quantum
algebras \cite{chak}. Let us recall certain notations of $(p,q)$-calculus:

\parindent8mm The $(p,q)$-integers $[n]_{p,q}$ are defined by

\begin{equation*}
[n]_{p,q}:=\frac{p^n-q^n}{p-q},~~~n=0,1,2,\dots,~~0<q<p\leq1.
\end{equation*}

The $(p,q)$-facorial and $(p,q)$-Binomial coefficients are defined by
\begin{equation*}
    [n]_{p,q}!:=\left\{
                  \begin{array}{ll}
                    [n]_{p,q}[n-1]_{p,q}\dots[1]_{p,q}, & \hbox{$n\in\mathbb{N}$;} \\
                    1, & \hbox{$n=0$}
                  \end{array}
                \right.
\end{equation*}
and
\begin{equation*}
\left[
\begin{array}{c}
n \\
k%
\end{array}%
\right] _{p,q}:=\frac{[n]_{p,q}!}{[k]_{p,q}![n-k]_{p,q}!},
\end{equation*}
respectively. Further, the $(p,q)$-binomial expansions are given as
\begin{equation*}
    (ax+by)_{p,q}^n:=\sum\limits_{k=0}^{n}p^{\binom{n-k}{2}}q^{\binom{k}{2}}a^{n-k}b^kx^{n-k}y^k
\end{equation*}
and
\begin{equation*}
(x-y)_{p,q}^n:=(x-y)(px-qy)(p^2x-q^2y)\cdots(p^{n-1}x-q^{n-1}y).
\end{equation*}
Let $m$ and $n$ be two non-negative integers. Then the following assertion is valid
\begin{equation*}
    (x-y)_{p,q}^{m+n}:=(x-y)_{p,q}^{m}(p^mx-q^my)_{p,q}^{n}.
\end{equation*}
Also, the $(p,q)$-derivative of a function $f$, denoted by $D_{p,q}f$, is defined by
\begin{equation*}
    (D_{p,q}f)(x):=\frac{f(px)-f(qx)}{(p-q)x},~~x\ne0,~~(D_{p,q}f)(0):=f^\prime(0)
\end{equation*}
provided that $f$ is differentiable at $0$.The $(p,q)$-derivative fulfils the following product rules
\begin{eqnarray*}
  D_{p,q}(f(x)g(x))&:=&f(px)D_{p,q}g(x)+g(qx)D_{p,q}f(x),\\
  D_{p,q}(f(x)g(x))&:=&f(px)D_{p,q}g(x)+g(qx)D_{p,q}f(x).
\end{eqnarray*}
Moreover,
\begin{eqnarray*}
  D_{p,q}\bigg(\frac{f(x)}{g(x)}\bigg)&:=&\frac{g(qx)D_{p,q}f(x)-f(qx)D_{p,q}g(x)}{g(px)g(qx)},\\
  D_{p,q}\bigg(\frac{f(x)}{g(x)}\bigg)&:=&\frac{g(px)D_{p,q}f(x)-f(px)D_{p,q}g(x)}{g(px)g(qx)}.
\end{eqnarray*}
We consider the $(p,q)$-exponential functions in the following forms:
\begin{eqnarray*}
  e_{p,q}(x) &=& \sum\limits_{n=0}^{\infty}p^{n(n-1)/2}\frac{x^n}{[n]_{p,q}!}, \\
  E_{p,q}(x) &=& \sum\limits_{n=0}^{\infty}q^{n(n-1)/2}\frac{x^n}{[n]_{p,q}!},
\end{eqnarray*}
which satisfy the equality $e_{p,q}(x)E_{p,q}(-x)=1$.
The definite integrals of the function $f$ are defined by
\begin{equation*}
\int_{0}^{a}f(x)d_{p,q}x=(q-p)a\sum\limits_{k=0}^{\infty}\frac{p^k}{q^{k+1}}%
f\left(\frac{p^k}{q^{k+1}}a\right),~~~\text{when}~~~\left|\frac pq\right|<1,
\end{equation*}
and
\begin{equation*}
\int_{0}^{a}f(x)d_{p,q}x=(p-q)a\sum\limits_{k=0}^{\infty}\frac{q^k}{p^{k+1}}%
f\left(\frac{q^k}{p^{k+1}}a\right),~~~\text{when}~~~\left|\frac pq\right|>1.
\end{equation*}

Details on $(p,q)$-calculus can be found in \cite{chak,jag,sad}.
For $p=1$, all the notions of $(p,q)$-calculus are reduced to $q$-calculus.

\vspace{.5cm}
\noindent\textbf{2. Operators and estimation of moments}\\
\\
Now we set the $(p,q)$-Sz\'{a}sz-Mirakjan basis function as
\begin{equation*}
    s_n(p,q;x)=:E_{p,q}\big(-[n]_{p,q}x\big)\sum\limits_{k=0}^{\infty}q^{\frac{k(k-1)}{2}}\frac{\big([n]_{p,q}x\big)^k}{[k]_{p,q}!}.
\end{equation*}
For $q\in(0,p)$, $p\in(q,1]$ and $x\in[0,\infty)$, $s_n(p,q;x)\ge0$. We can easily check that
\begin{equation*}
    \sum\limits_{k=0}^{\infty}s_n(p,q;x)=:E_{p,q}\big(-[n]_{p,q}x\big)
    \sum\limits_{k=0}^{\infty}q^{\frac{k(k-1)}{2}}\frac{\big([n]_{p,q}x\big)^k}{[k]_{p,q}!}=1.
\end{equation*}
For $0<q<p\le1$ the $(p,q)$-Sz\'{a}sz-Mirakjan operators are defined as
\begin{equation*}
S_n(f,p,q;x)=[n]_{p,q}\sum\limits_{k=0}^{\infty}p^{-k}q^ks_{n,k}(p,q;x)f\bigg(\frac{[k]_{p,q}}{q^{k-1}[n]_{p,q}}\bigg),~~x\in[0,\infty).\eqno(1)
\end{equation*}

\parindent8mm From the definition of the $(p,q)$-Sz\'{a}sz-Mirakjan operators we derive the following formulas.\\
\\
\noindent\textbf{Lemma 1.} Let $0<q<p\le1$. We have
\begin{itemize}
  \item [(i)] $S_n(1,p,q;x)=1$;
  \item [(ii)] $S_n(t,p,q;x)=x$;
  \item [(iii)] $S_n(t^2,p,q;x)=\frac{px^2}{q}+\frac{x}{[n]_{p,q}}$;
  \item [(iv)] $S_n(t^3,p,q;x)=\frac{p^3}{q^3}x^3+\frac{p^2+2pq}{q[n]_{p,q}}x^2+\frac{q^2}{[n]_{p,q}^2}x$;
  \item [(v)] $S_n(t^4,p,q;x)=\frac{p^6}{q^6}x^4+\frac{p^3(p^2+2q+3q^2)}{q^4[n]_{p,q}}x^3+\frac{p(p^2+3pq+3q^2)}{q[n]_{p,q}}x^2+\frac{q^3}{[n]_{p,q}^3}x$.
\end{itemize}
Now we propose our Kantorovich variant of $(p,q)$-Sz\'{a}sz-Mirakjan operators (1) as follows:

\parindent8mm For $f\in C[0,\infty)$, $0<q<p\le1$ and each positive integer $n$,
\begin{equation*}
    K_n(f,p,q;x)=[n]_{p,q}\sum\limits_{k=0}^{\infty}p^{-k}q^ks_{n,k}(p,q;x)
    \int_{\frac{[k]_{p,q}}{q^{k-1}[n]_{p,q}}}^{\frac{[k+1]_{p,q}}{q^{k}[n]_{p,q}}}f(t)d_{p,q}t\eqno(2)
\end{equation*}
We will derive the recurrence formula for $K_n(t^m,p,q;x)$ and calculate $K_n(t^m,p,q;x)$ for $m=0,1,2$.\\
\\
\noindent\textbf{Lemma 2.} For the operators $K_n$ we have
\begin{equation*}
    K_n(t^m,p,q;x)=\frac1{[m+1]_{p,q}}\sum\limits_{j=0}^{m}\sum\limits_{i=0}^{j}\frac{p^i}{q^i[n]_{p,q}^{j-i}}\binom {j}{i}S_n(t^{m+i-j},p,q;x).\eqno(3)
\end{equation*}

\noindent\textbf{Proof.} Using the expansion $a^{m+1}-b^{m+1}=(a-b)(a^m+a^{m-1}b+\cdots+ab^{m-1}+b^m)$ we have
\begin{eqnarray*}
  \int_{\frac{[k]_{p,q}}{q^{k-1}[n]_{p,q}}}^{\frac{[k+1]_{p,q}}{q^{k}[n]_{p,q}}}t^m d_{p,q}t&=& \frac1{[m+1]_{p,q}}\Bigg\{\bigg(\frac{[k+1]_{p,q}}{q^k[n]_{p,q}}\bigg)^{m+1}-\bigg(\frac{[k]_{p,q}}{q^{k-1}[n]_{p,q}}\bigg)^{m+1}\Bigg\}.
\end{eqnarray*}
Using $[k+1]_{p,q}=p^k+q[k]_{p,q}$ and also $[k+1]_{p,q}=q^k+p[k]_{p,q}$, we have
\begin{eqnarray*}
  \int_{\frac{[k]_{p,q}}{q^{k-1}[n]_{p,q}}}^{\frac{[k+1]_{p,q}}{q^{k}[n]_{p,q}}}t^m d_{p,q}t&=& \frac1{[m+1]_{p,q}}\frac{p^k}{q^k[n]_{p,q}}\sum\limits_{j=0}^{m}\bigg(\frac{[k+1]_{p,q}}{q^k[n]_{p,q}}\bigg)^{j}\bigg(\frac{[k]_{p,q}}{q^{k-1}[n]_{p,q}}\bigg)^{m-j}\\
  &=& \frac1{[m+1]_{p,q}}\frac{p^k}{q^k[n]_{p,q}}\sum\limits_{j=0}^{m}\bigg(\frac{q^k+p[k]_{p,q}}{q^k[n]_{p,q}}\bigg)^{j}\bigg(\frac{[k]_{p,q}}{q^{k-1}[n]_{p,q}}\bigg)^{m-j}\\
  &=& \frac1{[m+1]_{p,q}}\frac{p^k}{q^k[n]_{p,q}}\sum\limits_{j=0}^{m}\sum\limits_{i=0}^{j}\binom {j}{i}\frac{p^i[k]_{p,q}^iq^{k(j-i)}}{q^{kj}[n]_{p,q}^j}\frac{[k]_{p,q}^{m-j}}{q^{(k-1)(m-j)}[n]_{p,q}^{m-j}}\\
  &=& \frac1{[m+1]_{p,q}}\frac{p^k}{q^k[n]_{p,q}}\sum\limits_{j=0}^{m}\sum\limits_{i=0}^{j}\binom {j}{i}\frac{p^i[k]_{p,q}^{m+i-j}}{q^{ki}[n]^m_{p,q}q^{(k-1)(m-j)}}.
\end{eqnarray*}
Writing this in the definition of $K_n(t^m,p,q;x)$, we get
\begin{eqnarray*}
  K_n(t^m,p,q;x) &=& [n]_{p,q}\sum\limits_{k=0}^{\infty}p^{-k}q^ks_{n,k}(p,q;x)\int_{\frac{[k]_{p,q}}{q^{k-1}[n]_{p,q}}}^{\frac{[k+1]_{p,q}}{q^{k}[n]_{p,q}}}t^md_{p,q}t \\
   &=& \frac1{[m+1]_{p,q}}\sum\limits_{j=0}^{m}\sum\limits_{k=0}^{\infty}s_{n,k}(p,q;x)\sum\limits_{i=0}^{j}\frac{p^i}{q^i[n]_{p,q}^{j-i}}\binom {j}{i}\frac{[k]_{p,q}^{m+i-j}}{q^{(k-1)(m+i-j)}[n]^{m+i-j}_{p,q}}\\
   &=& \frac1{[m+1]_{p,q}}\sum\limits_{j=0}^{m}\sum\limits_{i=0}^{j}\frac{p^i}{q^i[n]_{p,q}^{j-i}}\binom {j}{i}\sum\limits_{k=0}^{\infty}\frac{[k]_{p,q}^{m+i-j}}{q^{(k-1)(m+i-j)}[n]^{m+i-j}_{p,q}}s_{n,k}(p,q;x)\\
   &=&\frac1{[m+1]_{p,q}}\sum\limits_{j=0}^{m}\sum\limits_{i=0}^{j}\frac{p^i}{q^i[n]_{p,q}^{j-i}}\binom {j}{i}S_n(t^{m+i-j},p,q;x).
\end{eqnarray*}

Using the recurrence formula (3) we may easily calculate $K_n(t^m,p,q;x)$ for $m=0,1,2$.\\

\vspace{.5cm}
\noindent\textbf{Lemma 3.} We have
\begin{itemize}
  \item [(i)] $K_n(1,p,q;x)=1$;
  \item [(ii)] $K_n(t,p,q;x)=\frac{1}{q}x+\frac{1}{[2]_{p,q}[n]_{p,q}}$;
  \item [(iii)] $K_n(t^2,p,q;x)=\frac{p}{q^3}x^2+\Big(\frac{p+[2]_{p,q}}{q[3]_{p,q}[n]_{p,q}}+\frac{1}{q^2[n]_{p,q}}\Big)x+\frac1{[3]_{p,q}[n]^2_{p,q}}$;
  \item [(iv)] $K_n(t^3,p,q;x)=\frac{p^3}{q^6}x^3+\Big(\frac{p^2+2pq}{q^4[n]_{p,q}}+\frac{p(3p^2+2pq+q^2)}{q^3[4]_{p,q}[n]_{p,q}}\Big)x^2
  +\Big(\frac{1}{q[n]^2_{p,q}}+\frac{3p^2+2pq+q^2}{q^2[4]_{p,q}[n]^2_{p,q}}+\frac{3p+q}{q[4]_{p,q}[n]^2_{p,q}}\Big)x+\frac{1}{[4]_{p,q}[n]^3_{p,q}}$;
  \item [(v)] $K_n(t^4,p,q;x)=\frac{p^6}{q^{10}}x^4+\Big(\frac{p^3(p^2+2q+3q^2)}{q^8[n]_{p,q}}+\frac{p^3(4p^3+3p^2q+2pq^2+q^3)}{q^6[5]_{p,q}[n]_{p,q}}\Big)x^3
  +\Big(\frac{p(p^2+3pq+3q^2)}{q^5[n]^2_{p,q}}+\frac{(p^2+2pq)(4p^3+3p^2q+2pq^2+q^3)}{q^4[5]_{p,q}[n]^2_{p,q}}
  +\frac{p(6p^2+3pq+q^2)}{q^3[5]_{p,q}[n]^2_{p,q}}\Big)x^2+\Big(\frac1{q[n]^3_{p,q}}+\frac{4p^3+3p^2q+2pq^2+q^3}{q[5]_{p,q}[n]^3_{p,q}}
  +\frac{6p^2+3pq+q^2}{q^2[5]_{p,q}[n]^3_{p,q}}+\frac{4p+q}{q[5]_{p,q}[n]^3_{p,q}}\Big)x
  +\frac{1}{[5]_{p,q}[n]^4_{p,q}}$;
  \item [(iv)] $K_n\big((t-x),p,q;x\big)=\frac{1-q}{q}x+\frac{1}{[2]_{p,q}[n]_{p,q}}$;
  \item [(v)] $K_n\big((t-x)^2,p,q;x\big)=\Big(\frac{p}{q^3}-\frac2q+1\Big)x^2+\Big(\frac{p+[2]_{p,q}}{q[3]_{p,q}[n]_{p,q}}
      +\frac{1}{q^2[n]_{p,q}}-\frac{2}{[2]_{p,q}[n]_{p,q}}\Big)x+\frac1{[3]_{p,q}[n]^2_{p,q}};\hspace{2.5cm}(4)$
  \item [(vi)] $K_n\big((t-x)^4,p,q;x\big)=x^4\Big(\frac{p^6}{q^{10}}-\frac{4p^3}{q^6}+\frac{6p}{q^3}-\frac{4}{q}+1\Big)
  +\frac{x^3}{[n]_{p,q}}\Big(\frac{p^3(p^2+2q+3q^2)}{q^8}+\frac{p^3(4p^3+3p^2q+2pq^2+q^3)}{q^6[5]_{p,q}}
  -\frac{4(p^2+2pq)}{q^4}-\frac{4p(3p^2+2pq+q^2)}{q^3[4]_{p,q}}+\frac{6(2p+q)}{q[3]_{p,q}}+\frac6{q^2}-\frac4{[2]_{p,q}}\Big)
  +\frac{x^2}{[n]^2_{p,q}}\Big(\frac{p(p^2+3pq+3q^2)}{q^5}+\frac{(p^2+2pq)(4p^3+3p^2q+2pq^2+q^3)}{q^4[5]_{p,q}}
  +\frac{p(6p^2+3pq+q^2)}{q^3[5]_{p,q}}-\frac{4(3p^2+2pq+q^2)}{q^2}-\frac4{q}+\frac6{[3]_{p,q}}\Big)
  +\frac{x}{[n]^3_{p,q}}\Big(\frac{1}{q}+\frac{4p^3+3p^2q+2pq^2+q^3}{q[5]_{p,q}}
  +\frac{6p^2+3pq+q^2}{q^2[5]_{p,q}}+\frac{4p+q}{q[5]_{p,q}}\Big)$.
\end{itemize}

\noindent\textbf{Proof.} Obviously, with the help of Lemma 1, we can get

\begin{eqnarray*}
  K_n(t,p,q;x) &=& \frac{1}{[2]_{p,q}}\bigg\{\Big(1+\frac pq\Big)S_n(t,p,q;x)+\frac1{[n]_{p,q}}S_n(1,p,q;x)\bigg\} \\
  &=& \frac{1}{q}x+\frac{1}{[2]_{p,q}[n]_{p,q}},
\end{eqnarray*}

\begin{eqnarray*}
  K_n(t^2,p,q;x) &=& \frac{1}{[3]_{p,q}}\bigg\{\Big(1+\frac pq+\frac{p^2}{q^2}\Big)S_n(t^2,p,q;x)+\Big(\frac1{[n]_{p,q}}+\frac {2p}{q[n]_{p,q}}\Big)S_n(t,p,q;x)+\frac1{[n]^2_{p,q}}S_n(1,p,q;x)\bigg\} \\
  &=& \frac{1}{q^2}S_n(t^2,p,q;x)+\frac{p+[2]_{p,q}}{q[3]_{p,q}[n]_{p,q}}S_n(t,p,q;x)+\frac1{[3]_{p,q}[n]^2_{p,q}}S_n(1,p,q;x)\\
  &=& \frac{p}{q^3}x^2+\Big(\frac{p+[2]_{p,q}}{q[3]_{p,q}[n]_{p,q}}+\frac{1}{q^2[n]_{p,q}}\Big)x+\frac1{[3]_{p,q}[n]^2_{p,q}}.
\end{eqnarray*}

\vspace{.5cm}
Using the linearity of the operators, we can have
\begin{eqnarray*}
  K_n\big((t-x)^2,p,q;x\big) &=& K_n(t^2,p,q;x)-2xK_n(t,p,q;x)+x^2K_n(1,p,q;x)\\
  &=& \frac{p}{q^3}x^2+\Big(\frac{p+[2]_{p,q}}{q[3]_{p,q}[n]_{p,q}}+\frac{1}{q^2[n]_{p,q}}\Big)x+\frac1{[3]_{p,q}[n]^2_{p,q}}
  -2x\Big(\frac{1}{q}x+\frac{1}{[2]_{p,q}[n]_{p,q}}\Big)+x^2\\
  &=&\Big(\frac{p}{q^3}-\frac2q+1\Big)x^2+\Big(\frac{p+[2]_{p,q}}{q[3]_{p,q}[n]_{p,q}}
      +\frac{1}{q^2[n]_{p,q}}-\frac{2}{[2]_{p,q}[n]_{p,q}}\Big)x+\frac1{[3]_{p,q}[n]^2_{p,q}}.
\end{eqnarray*}

\vspace{.5cm}
\noindent\textbf{Remark} For $q\in(0,1)$ and $p\in(q,1]$ it is obvious that
$\lim\limits_{n\to\infty}[n]_{p,q}=\frac1{p-q}$. In order to reach to convergence
results of the operator $K_n$ we take sequences $q_n\in(0,1)$ and $p_n\in(q_n,1]$
such that $\lim\limits_{n\to\infty}p_n=1$ $\lim\limits_{n\to\infty}q_n=1$. So we get that
$\lim\limits_{n\to\infty}[n]_{p_n,q_n}=\infty$.\newline

Thus the above remark provides an example that such a sequence can always be constructed.
If we choose for $a>b>0$, $q_n=\frac{n}{n+a}<\frac{n}{n+b}=p_n$
such that $0<q_n<p_n\leq1$, it can be easily seen that
 $\lim\limits_{n\to\infty}p_n=1,~\lim\limits_{n\to\infty}q_n=1$ and
$\lim\limits_{n\to\infty}p_n^n=e^{-b},~\lim\limits_{n\to\infty}q_n^n=e^{-a}$.
Hence we guarantee that $\lim\limits_{n\to\infty}[n]_{p_n,q_n}=\infty$.


\vspace{.5cm}
\noindent\textbf{3. Direct approximation result}

\vspace{.5cm}
In this section we study the Korovkin's approximation property of the Kantorovich variant of $(p,q)$-Sz\'{a}sz operators.

\vspace{.5cm}
\noindent\textbf{Theorem 4.} Let $0<q_n<p_n\le1$ and $A>0$. Then for each $f\in C_m[0,\infty):=\big\{f\in C[0,\infty):|f(x)|\le M_f(1+x^m),~\text{for some}~M_f>0~\text{depending on}~f,~m>0\big\}$ where $C_m[0,\infty)$ be endowed with the norm $\|f\|_m=\sup\limits_{x\in[0,\infty)}\frac{|f(x)|}{1+x^m}$, the sequence of operators $K_n(f,p_n,q_n;x)$ converges to $f$ uniformly on $[0,A]$ if and only if $\lim_{n\to\infty}p_n=1$ and $\lim_{n\to\infty}q_n=1$.

\vspace{.5cm}
\noindent\textbf{Proof.} First, we assume that $\lim_{n\to\infty}p_n=1$ and $\lim_{n\to\infty}q_n=1$. Now, we have to show that $K_n(f,p_n,q_n;x)$ converges to $f$ uniformly on $[0,A]$.\\

\parindent8mm From Lemma 3, we see that
\begin{equation*}
    K_n(1,p_n,q_n;x)\to1,~~~K_n(t,p_n,q_n;x)\to x,~~~K_n(t^2,p_n,q_n;x)\to x^2,
\end{equation*}

uniformly on $[0,A]$ as $n\to\infty$.\\

\parindent8mm Therefore, the well-known property of the Korovkin theorem implies that $K_n(f,p_n,q_n;x)$ converges to $f$ uniformly on $[0,A]$ provided $f\in C_m[0,\infty)$.\\

\parindent8mm We show the converse part by contradiction. Assume that $p_n$ and $q_n$ do not converge to 1. Then they must contain subsequences $p_{n_k}\in (0,1)$, $q_{n_k}\in (0,1)$, $p_{n_k}\to a\in [0,1)$ and $q_{n_k}\to b\in [0,1)$ as $k\to\infty$, respectively.\\

Thus,
\begin{equation*}
    \frac{1}{[n_k]_{p_{n_k},q_{n_k}}}=\frac{p_{n_k}-q_{n_k}}{(p_{n_k})^{n_k}-(q_{n_k})^{n_k}}\to 0 \text{ as } k\to\infty
\end{equation*}
and we get

\begin{equation*}
    K_n(t,p_{n_k},q_{n_k};x)-x=\frac{1}{q_{n_k}}x+\frac{1}{[2]_{p_{n_k},q_{n_k}}[n]_{p_{n_k},q_{n_k}}}-x\to \frac{x}{b}-x\ne0.
\end{equation*}
This leads to a contradiction. Thus $p_n\to1$ and $q_n\to1$ as $n\to\infty$.\\

\vspace{.5cm}
\noindent\textbf{Theorem 5.} Let $f\in C_{2}[0,\infty),~q=q_n\in(0,1)$ and $p=p_n\in(q,1]$ such that $p_n\to1,~q_n\to1$ as $n\to\infty$ and $\omega_{a+1}(f,\delta)$ be the modulus of continuity on the finite interval
$[0,a+1]\subset[0,\infty)$, where $a>0$. Then
\begin{equation*}
\left|K_n(f,p,q;x)-f(x)\right|\leq 4M_f(1+a^2)\delta^2_n(x)+2\omega_{a+1}(f,\delta_n(x))
\end{equation*}
where $\delta_n(x)=\sqrt{K_n{\big((t-x)^2,p,q;x\big)}}$, given by (4).\\

\noindent\textbf{Proof.} For $x\in[0,a]$ and $t>a+1$, since $t-x>1$, we have
\begin{equation*}
    |f(t)-f(x)|\leq M_f(2+x^2+t^2)\leq M_f\bigl{(}2+3x^2+2(t-x)^2\bigl{)}\leq M_f(4+3x^2)(t-x)^2\leq4M_f(1+a^2)(t-x)^2.\eqno(5)
\end{equation*}
For $x\in[0,a]$ and $t\leq a+1$, we have
\begin{equation*}
    |f(t)-f(x)|\leq\omega_{a+1}(f,|t-x|)\leq\left(1+\frac{|t-x|}{\delta}\right)\omega_{a+1}(f,\delta)\eqno(6)
\end{equation*}
with $\delta>0$.\\

\parindent8mm From (5) and (6), we may write

\begin{equation*}
    |f(t)-f(x)|\leq4M_f(1+a^2)(t-x)^2+\left(1+\frac{|t-x|}{\delta}\right)\omega_{a+1}(f,\delta), \eqno(7)
\end{equation*}
for $x\in[0,a]$ and $t\geq0$. Thus, by applying the Cauchy-Schwarz's inequality, we have
\begin{eqnarray*}
  \left|K_n(f,p,q;x)-f(x)\right|&\leq& K_n\bigl{(}|f(t)-f(x)|,p,q;x\bigl{)}\\
  &\leq&4M_f(1+a^2) K_n\bigl{(}(t-x)^2,p,q;x\bigl{)}
  +\left(1+\frac1\delta\sqrt{K_n\bigl{(}(t-x)^2,p,q;x\bigl{)}}\right)\omega_{a+1}(f,\delta)\\
   &\le&4M_f(1+a^2)\delta^2_n(x)+2\omega_{a+1}(f,\delta_n(x))
\end{eqnarray*}
  on choosing $\delta:=\delta_n(x)$. This completes the proof of the theorem.\\

\vspace{.5cm}
\noindent\textbf{4. Local approximation}

\vspace{.5cm}
In this section we establish local approximation theorem for the Kantorovich variant of $(p,q)$-Sz\'{a}sz operators. Let $C_B[0,\infty)$ be the space of all real-valued continuous bounded functions $f$ on $[0,\infty)$, endowed with the norm $\|f\|=\sup\limits_{x\in[0,\infty)}|f(x)|$. The Peetre's K-functional is defined by

\begin{equation*}
    K_2(f,\delta)=\inf\limits_{g\in C^2[0,\infty)}\big\{\|f-g\|+\delta\|g^{\prime\prime}\|\big\},
\end{equation*}
where $C_B^2[0,\infty):=\big\{g\in C_B[0,\infty):g^\prime,g^{\prime\prime}\in C_B[0,\infty)\big\}$. By [4, p.177, Theorem 2.4], there exists an absolute constant $M>0$ such that

\begin{equation*}
    K_2(f,\delta)\leq M\omega_2(f,\sqrt{\delta}),\eqno(8)
\end{equation*}
where $\delta>0$ and the second order modulus of smoothness is defined as

\begin{equation*}
    \omega_2(f,\delta)=\sup\limits_{0<h\le\delta}\sup\limits_{x\in[0,\infty)}|f(x+2h)-2f(x+h)+f(x)|,\eqno(9)
\end{equation*}
where $f\in C_B[0,\infty)$ and $\delta>0$.\\

\vspace{.5cm}
\noindent\textbf{Theorem 6.} Let $f\in C_B[0,\infty)$ and $0<q<p\le1$. Then, for every $x\in[0,\infty)$, we have

\begin{equation*}
    |K_n(f,p,q;x)-f(x)|\leq M\omega_2\Big(f,\sqrt{\delta_n(x)}\Big)+\omega\Big(f,\frac{1}{[2]_{p,q}[n]_{p,q}}+\frac{1-q}{q}x\Big),
\end{equation*}
where $M$ is an absolute constant and

\begin{equation*}
    \delta_n(x)=K_n\big((t-x)^2,p,q;x\big)+\Big(\frac{1}{[2]_{p,q}[n]_{p,q}}+\frac{1-q}{q}x\Big)^2.
\end{equation*}

\vspace{.5cm}
\noindent\textbf{Proof.} For $x\in[0,\infty)$, we consider the auxiliary operators $K_n^*$ defined by
\begin{equation*}
    K_n^*(f,p,q;x)=K_n(f,p,q;x)-f\Big(\frac{1}{[2]_{p,q}[n]_{p,q}}+\frac{1}{q}x\Big)+f(x).\eqno(10)
\end{equation*}
From Lemma 3, we observe that the operators $K_n^*(f,p,q;x)$ are linear and reproduce the linear functions.
Hence
\begin{eqnarray*}
  K_n^*\big((t-x),p,q;x\big) &=& K_n\big((t-x),p,q;x\big)-\Big(\frac{1}{[2]_{p,q}[n]_{p,q}}+\frac{1}{q}x-x\Big) \\
   &=& K_n(t,p,q;x)-xK_n(1,p,q;x)-\Big(\frac{1}{[2]_{p,q}[n]_{p,q}}+\frac{1}{q}x\Big)+x=0.\hspace{1.5cm}(11)
\end{eqnarray*}
Let $x\in[0,\infty)$ and $g\in C_B^2[0,\infty)$. Using the Taylor's formula
\begin{equation*}
g(t)=g(x)+g^{\prime }(x)(t-x)+\int_{x}^{t}(t-u)~g^{\prime \prime
}(u)~du.
\end{equation*}%
\newline
Applying $K_n^*$ to both sides of the above equation and using (11), we have
\begin{eqnarray*}
  K_n^*(g,p,q;x)-g(x) &=& K_n^*\big((t-x)g^\prime(x),p,q;x\big)+K_n^*\bigg(\int_x^t(t-u)g^{\prime\prime}(u)du,p,q;x\bigg) \\
   &=& g^\prime(x)K_n^*\big((t-x),p,q;x\big)+K_n^{(p,q)}\bigg(\int_x^t(t-u)g^{\prime\prime}(u)du,p,q;x\bigg)\\
   &&-\int_x^{\frac{1}{[2]_{p,q}[n]_{p,q}}+\frac{1}{q}x}
   \Big(\frac{1}{[2]_{p,q}[n]_{p,q}}+\frac{1}{q}x-u\Big)g^{\prime\prime}(u)du\\
   &=&K_n\bigg(\int_x^t(t-u)g^{\prime\prime}(u)du,p,q;x\bigg)\\
   &&-\int_x^{\frac{1}{[2]_{p,q}[n]_{p,q}}+\frac{1}{q}x}
   \Big(\frac{1}{[2]_{p,q}[n]_{p,q}}+\frac{1}{q}x-u\Big)g^{\prime\prime}(u)du.
\end{eqnarray*}
On the other hand, since
\begin{equation*}
    \bigg|\int_x^t(t-u)g^{\prime\prime}(u)du\bigg|\leq\int_x^t|t-u||g^{\prime\prime}(u)|du\leq\|g^{\prime\prime}\|\int_x^t|t-u|du\leq(t-x)^2\|g^{\prime\prime}\|
\end{equation*}
and
\begin{eqnarray*}
  \Bigg|\int_x^{\frac{1}{[2]_{p,q}[n]_{p,q}}+\frac{1}{q}x}
   \Big(\frac{1}{[2]_{p,q}[n]_{p,q}}+\frac{1}{q}x-u\Big)g^{\prime\prime}(u)du\Bigg|\\
   \leq \Big(\frac{1}{[2]_{p,q}[n]_{p,q}}+\frac{1}{q}x-u\Big)^2\|g^{\prime\prime}\|  \\
\end{eqnarray*}
we conclude that
\begin{eqnarray*}
  \Big|K_n^*(g,p,q;x)-g(x)\Big| &=& \bigg|K_n\bigg(\int_x^t(t-u)g^{\prime\prime}(u)du,p,q;x\bigg)\\
   &&-\int_x^{\frac{1}{[2]_{p,q}[n]_{p,q}}+\frac{1}{q}x}
   \Big(\frac{1}{[2]_{p,q}[n]_{p,q}}+\frac{1}{q}x-u\Big)g^{\prime\prime}(u)du\bigg| \\
   &\leq& \Vert g^{\prime \prime }\Vert K_{n}\big((t-x)^{2},p,q;x\big)+\Big(\frac{1}{[2]_{p,q}[n]_{p,q}}+\frac{1}{q}x-x\Big)^2\Vert g^{\prime \prime }\Vert\\
   &=&\delta_n(x)\Vert g^{\prime \prime }\Vert.\hspace{9cm}(12)
\end{eqnarray*}
Now, taking into account boundedness of $K_n^*$ by (10), we have
\begin{equation*}
    \big|K_n^*(f,p,q;x)\big|\leq\big|K_n(f,p,q;x)\big|+2\|f\|\leq3\|f\|\eqno(13)
\end{equation*}
Using (12) and (13) in (10), we obtain
\begin{eqnarray*}
  \big|K_n(f,p,q;x)-f(x)\big| &\leq& \big|K_n^*(f,p,q;x)-f(x)\big|+\bigg|f(x)-f\Big(\frac{1}{[2]_{p,q}[n]_{p,q}}+\frac{1}{q}x\Big)\bigg| \\
   &\leq& \big|K_n^*(f-g,p,q;x)-(f-g)(x)\big|\\
   &&+\bigg|f(x)-f\Big(\frac{1}{[2]_{p,q}[n]_{p,q}}+\frac{1}{q}x\Big)\bigg|
   +\big|K_n^*(g,p,q;x)-g(x)\big|\\
   &\leq& \big|K_n^*(f-g,p,q;x)\big|+\big|(f-g)(x)\big|\\
   &&+\bigg|f(x)-f\Big(\frac{1}{[2]_{p,q}[n]_{p,q}}+\frac{1}{q}x\Big)\bigg|
   +\big|K_n^*(g,p,q;x)-g(x)\big|\\
   &\leq&4\|f-g\|+\omega\Big(f,\frac{1}{[2]_{p,q}[n]_{p,q}}+\frac{1-q}{q}x\Big)
   +\delta_n(x)\Vert g^{\prime \prime }\Vert.
\end{eqnarray*}
Hence, taking the infimum on the right-hand side over all $g\in C_B^2[0,\infty)$, we have the following result
\begin{equation*}
    \big|K_n(f,p,q;x)-f(x)\big|\leq4K_2\big(f,\delta_n(x)\big)+\omega\Big(f,\frac{1}{[2]_{p,q}[n]_{p,q}}+\frac{1-q}{q}x\Big).
\end{equation*}
In view of the property of $K$-functional (8), we get
\begin{equation*}
    \big|K_n(f,p,q;x)-f(x)\big|\leq M\omega_2\Big(f,\sqrt{\delta_n(x)}\Big)+\omega\Big(f,\frac{1}{[2]_{p,q}[n]_{p,q}}+\frac{1-q}{q}x\Big).
\end{equation*}

\parindent=0mm This completes the proof of the theorem.\\

\vspace{.5cm}
\noindent\textbf{5. Weighted approximation}

\vspace{.5cm}
Let $f\in C_2^*[0,\infty):=\big\{f\in C_2[0,\infty):\lim\limits_{x\to\infty}\frac{|f(x)|}{1+x^2}<\infty\big\}$.
Throughout the section, we assume that $(p_n)$ and $(q_n)$ are sequences such that $0<q_n<p_n\le1$ and $p_n\to1$, $q_n\to1$ as $n\to\infty$.

\vspace{.5cm}
\textbf{Theorem 7.} For each $f\in C^*_{2}[0,\infty)$, we have
\begin{equation*}
    \lim\limits_{n\to\infty}\|K_n(f,p_n,q_n;x)-f(x)\|_{2}=0.
\end{equation*}
\\
\textbf{Proof.} Using the Korovkin type theorem on weighted approximation in \cite{gad} we see that it is sufficient to verify the following three conditions
\begin{equation*}
   \lim\limits_{n\to\infty}\|K_{n}(t^i,p_n,q_n;x)-x^i\|_{2}=0, ~~~i=0,1,2.\eqno(14)
\end{equation*}
Since $K_{n}(1,p_n,q_n;x)=1$, (14) holds true for $m=0$.\\

\parindent8mm By Lemma 3, we have
\begin{eqnarray*}
  \|K_{n}(t,p_n,q_n;x)-x\|_{2}&=&\sup\limits_{x\in[0,\infty)}\frac{|K_{n}(t,p_n,q_n;x)-x|}{1+x^2}\\
  &=&\sup\limits_{x\in[0,\infty)}\frac1{1+x^2}\left|\frac1{q_n}x+\frac1{[2]_{p,q}[n]_{p,q}}-x\right|\\
  &\leq&\Big(\frac{1}{q_n}-1\Big)\sup\limits_{x\in[0,\infty)}\frac{x}{1+x^2}+\frac1{[2]_{p,q}[n]_{p,q}}\sup\limits_{x\in[0,\infty)}\frac{1}{1+x^2}\\
  &\leq&\frac{1}{q_n}-1+\frac1{[2]_{p,q}[n]_{p,q}}.
\end{eqnarray*}
which implies that the condition in (14) holds for $i=1$ as $n\to\infty$.\\

\noindent Similarly we can write

\begin{eqnarray*}
  \|K_{n}(t^2,p_n,q_n;x)-x^2\|_{2}&=&\sup\limits_{x\in[0,\infty)}\frac{|K_{n}(t^2,p_n,q_n;x)-x^2|}{1+x^2}\\
  &\leq&\Big(\frac{p_n}{q_n^3}-1\Big)\sup\limits_{x\in[0,\infty)}\frac{x^2}{1+x^2}\\
  &&+\Big(\frac{2p_n+q_n}{q_n{[3]_{p_n,q_n}}{[n]_{p_n,q_n}}}+\frac{1}{q_n^2[n]_{p_n,q_n}}\Big)\sup\limits_{x\in[0,\infty)}\frac{x}{1+x^2}\\
  &&+\frac{1}{[3]_{p_n,q_n}[n]^2_{p_n,q_n}}\sup\limits_{x\in[0,\infty)}\frac{1}{1+x^2}\\
  &\leq&\frac{p_n}{q_n^3}-1+\frac{2p_n+q_n}{q_n{[3]_{p_n,q_n}}{[n]_{p_n,q_n}}}+\frac{1}{q_n^2[n]_{p_n,q_n}}+\frac{1}{[3]_{p_n,q_n}[n]^2_{p_n,q_n}},
\end{eqnarray*}
which implies that
\begin{equation*}
    \lim\limits_{n\to\infty}\|K_{n}(t^2,p_n,q_n;x)-x^2\|_{2}=0,
\end{equation*}
(14) holds for $i=2$.
Thus the proof is completed.\\

\parindent8mm We give the following theorem to approximate all functions in $C^*_{2}[0,\infty)$. This type of results are given in \cite{dog}
for classical Sz\'{a}sz operators.\\
\\
\textbf{Theorem 8.} For each $f\in C^*_{2}[0,\infty)$ and $\alpha>0$, we have
\begin{equation*}
       \lim\limits_{n\to\infty}\sup\limits_{x\in[0,\infty)}\frac{\bigl{|}K_{n}(f,p_n,q_n;x)-f(x)\bigl{|}}{(1+x^2)^{1+\alpha}}=0.
\end{equation*}\\
\\
\textbf{Proof.} Let $x_0\in[0,\infty)$ be arbitrary but fixed. Then
\begin{align*}
    \sup\limits_{x\in[0,\infty)}\frac{\bigl{|}K_{n}(f,p_n,q_n;x)-f(x)\bigl{|}}{(1+x^2)^{1+\alpha}}
    =\sup\limits_{x\leq x_0}\frac{\bigl{|}K_{n}(f,p_n,q_n;x)-f(x)\bigl{|}}{(1+x^2)^{1+\alpha}}
    +\sup\limits_{x> x_0}\frac{\bigl{|}K_{n}(f,p_n,q_n;x)-f(x)\bigl{|}}{(1+x^2)^{1+\alpha}}\\
    \leq\|K_n(f)-f\|_{C[0,x_0]}+\|f\|_{2}\sup\limits_{x> x_0}\frac{\bigl{|}K_n\big(1+t^2,p,q;x\big)\bigl{|}}
    {(1+x^2)^{1+\alpha}}+\sup\limits_{x>x_0}\frac{|f(x)|}{(1+x^2)^{1+\alpha}}.~~(15)
\end{align*}
Since $|f(x)|\leq\|f\|_2(1+x^2)$, we have $\sup\limits_{x>x_0}\frac{|f(x)|}{(1+x^2)^{1+\alpha}}\le\frac{\|f\|_2}{(1+x_0^2)^\alpha}$.\\

\parindent8mm Let $\varepsilon>0$ be arbitrary. We can choose $x_0$ to be so large that
\begin{equation*}
    \frac{\|f\|_2}{(1+x_0^2)^\alpha}<\frac{\varepsilon}{3}.\eqno(16)
\end{equation*}
In view of Theorem 4, we obtain
\begin{equation*}
    \|f\|_{2}\lim\limits_{n\to\infty}\frac{\bigl{|}K_n\big(1+t^2,p,q;x\big)\bigl{|}}
    {(1+x^2)^{1+\alpha}}=\frac{1+x^2}{(1+x^2)^{1+\alpha}}\|f\|_2=\frac{\|f\|_2}{(1+x^2)^{\alpha}}
    \le\frac{\|f\|_2}{(1+x_0^2)^{\alpha}}<\frac{\varepsilon}{3}.\eqno(17)
\end{equation*}

Using Theorem 5, we can see that the first term of the inequality (15), implies that
\begin{equation*}
    \|K_n(f)-f\|_{C[0,x_0]}<\frac{\varepsilon}{3},~~~\text{as}~n\to\infty.\eqno(18)
\end{equation*}
Combining (16)-(18), we get that desired result.

\vspace{.5cm}
For $f\in C_2^*[0,\infty)$, the weighted modulus of continuity is defined as
\begin{equation*}
    \Omega_2(f,\delta)=\sup\limits_{x\ge0,0<h\le\delta}\frac{|f(x+h)-f(x)|}{1+(x+h)^2}.
\end{equation*}

\vspace{.5cm}
\noindent\textbf{Lemma 9} (\cite{lop}){\bf{.}} If $f\in C_2^*[0,\infty)$ then
\begin{itemize}
  \item [(i)] $\Omega_2(f,\delta)$ is monotone increasing function of $\delta$,
  \item [(ii)] $\lim\limits_{\delta\to0^+}\Omega_2(f,\delta)=0$,
  \item [(iii)] for any $\lambda\in[0,\infty),~\Omega_2(f,\lambda\delta)\leq(1+\lambda)\Omega_2(f,\delta).$
\end{itemize}

\vspace{.5cm}
\noindent\textbf{Theorem 10.} If $f\in C_2^*[0,\infty)$, then for sufficiently large $n$ we have
\begin{equation*}
    \big|K_n(f,p,q;x)-f(x)\big|\leq K(1+x^{2+\lambda})\Omega_2(f,\delta_n),~~x\in[0,\infty),
\end{equation*}
where $\lambda\ge1$, $\delta_n=\max\{\alpha_n,\beta_n,\gamma_n\}$, $\alpha_n,\beta_n,\gamma_n$ being
\begin{equation*}
    \alpha_n=\frac{p}{q^3}-\frac2q+1,~~~\beta_n=\frac{p+[2]_{p,q}}{q[3]_{p,q}[n]_{p,q}}
      +\frac{1}{q^2[n]_{p,q}}-\frac{2}{[2]_{p,q}[n]_{p,q}},~~~\gamma_n=\frac1{[3]_{p,q}[n]^2_{p,q}}.
\end{equation*}

\vspace{.5cm}
\noindent\textbf{Proof.} From the definition of $\Omega _{2}(f,\delta )$ and Lemma
9, we may write
\begin{eqnarray*}
|f(t)-f(x)| &\leq &\bigl{(}1+(x+|t-x|)^{m}\bigl{)}\left( \frac{|t-x|}{\delta
}+1\right) \Omega _{2}(f,\delta ) \\
&\leq &\bigl{(}1+(2x+t)^{m}\bigl{)}\left( \frac{|t-x|}{\delta }+1\right)
\Omega _{2}(f,\delta ) \\
:= &&\varphi_{x}(t)\left( 1+\frac{1}{\delta}\psi_x(t)\right) \Omega _{2}(f,\delta ).
\end{eqnarray*}%
\newline
Then we obtain
\begin{equation*}
|K_{n}(f,p,q;x)-f(x)|\leq \Omega _{2}(f,\delta_n )\left( K_{n}(\varphi _{x},p,q;x)+\frac1{\delta_n}K_{n}\left( \varphi_{x}\psi_x,p,q%
;x\right) \right) .
\end{equation*}%
Applying the Cauchy-Schwartz inequality to the second term on the right-hand side, we get
\begin{equation*}
|K_{n}(f,p,q;x)-f(x)|\leq \Omega _{2}(f,\delta )\Big( K_{n}(\varphi _{x},p,q;x)
+\frac1{\delta_n}\sqrt{K_{n}(\varphi _{x}^2,p,q;x)}\sqrt{K_{n}(\psi _{x}^2,p,q;x)}\Big) .\eqno(19)
\end{equation*}%
From Lemma 3, we get
\begin{eqnarray*}
  \frac1{1+x^2}K_{n}(1+t^2,p,q;x) &=& \frac1{1+x^2}+\frac{p}{q^3}\frac{x^2}{1+x^2}+\Big(\frac{p+[2]_{p,q}}{q[3]_{p,q}[n]_{p,q}}
      +\frac{1}{q^2[n]_{p,q}}\Big)\frac{x}{1+x^2}+\frac1{[3]_{p,q}[n]^2_{p,q}} \frac1{1+x^2}\\
   &\leq& 1+C_1,~\text{for sufficiently large $n$}\hspace{6cm}(20)
\end{eqnarray*}
where $C_1$ is a positive constant. From (20), there exists a positive constant $K_1$ such that $K_{n}(\varphi_x,p,q;x)\le K_1(1+x^2)$, for sufficiently large $n$.\\
Proceeding similarly $\frac1{1+x^4}K_{n}(1+t^4,p,q;x)\le 1+C_2$, for sufficiently large $n$, where $C_2$ is a positive constant.\\
So there exists a positive constant $K_2$ such that $K_{n}(\varphi_x^2,p,q;x)\le K_2(1+x^2)$, where $x\in[0,\infty)$ $n$ is large enough. Also we get
\begin{eqnarray*}
  K_n(\psi_x^2,p,q;x) &=& \Big(\frac{p}{q^3}-\frac2q+1\Big)x^2+\Big(\frac{p+[2]_{p,q}}{q[3]_{p,q}[n]_{p,q}}
      +\frac{1}{q^2[n]_{p,q}}-\frac{2}{[2]_{p,q}[n]_{p,q}}\Big)x+\frac1{[3]_{p,q}[n]^2_{p,q}}\\
  &=&\alpha_nx^2+\beta_nx+\gamma_n.
\end{eqnarray*}
Hence form (19), we have
\begin{equation*}
    |K_{n}(f,p,q;x)-f(x)|\leq (1+x^2)\Big(K_1+\frac1{\delta_n}K_2\sqrt{\alpha_nx^2+\beta_nx+\gamma_n}\Big)\Omega _{2}(f,\delta_n).
\end{equation*}
If we take $\delta_n=\max\{\alpha_n,\beta_n,\gamma_n\}$, then we get
\begin{eqnarray*}
     |K_{n}(f,p,q;x)-f(x)|&\leq& (1+x^2)\Big(K_1+K_2\sqrt{x^2+x+1}\Big)\Omega _{2}(f,\delta_n)\\
 &\leq& K_3(1+x^{2+\lambda})\Omega _{2}(f,\delta_n),~~~\text{for sufficiently large $n$ and $x\in[0,\infty)$}.
\end{eqnarray*}
Hence the proof is completed.

\end{document}